\documentclass[amssymb,amsfonts,refcheck,11pt,verbatim,righttag]{amsart}
\usepackage{amssymb}

\usepackage{graphicx}

\numberwithin{equation}{section}
\theoremstyle{plain}

\newtheorem{thm}{Theorem}[section]

\newtheorem{cor}[thm]{Corollary}
\newtheorem{ex}[thm]{Example}

\def\R {{\Bbb R}}
\def\N {{\Bbb N}}

\def\M {{\mathcal M}}

\def\A {{\mathcal A}}
\def\B {{\mathcal B}}
\def\C {{\mathcal C}}

\def\bq{{\bf q}}

\begin{document}
\baselineskip 16pt
\title{Non-uniqueness of ergodic measures with full Hausdorff dimension on a Gatzouras-Lalley carpet}

\author{Julien Barral}
\address{LAGA (UMR 7539), D\'epartement de Math\'ematiques, Institut Galil\'ee, Universit\'e Paris 13, 99 avenue Jean-Baptiste Cl\'ement , 93430  Villetaneuse, France}
\email{barral@math.univ-paris13.fr}
\author{De-Jun Feng}
\address{
Department of Mathematics\\
The Chinese University of Hong Kong\\
Shatin,  Hong Kong\\
} \email{djfeng@math.cuhk.edu.hk}

\keywords{Ergodic invariant measures, non-conformal repellers, Hausdorff dimension}
\thanks {
2000 {\it Mathematics Subject Classification}: Primary 37D35, Secondary 37B10, 37A35, 28A78}

\date{}

\begin{abstract}
In this note, we show that on certain  Gatzouras-Lalley carpet,  there exist more than one ergodic measures with full Hausdorff dimension. This gives a negative answer to a conjecture of Gatzouras and Peres in \cite{GaPe97}. \

\end{abstract}

\maketitle
\section{Introduction}
The problem we are interested in is  the uniqueness of ergodic invariant measures on non-conformal repellers with full Hausdorff dimension (see  \cite{GaPe96, ChPe10} for a survey). For $C^{1+\alpha}$ conformal repellers, the existence and the uniqueness of an ergodic measure with full dimension follows from Bowen's equation together with the classical thermodynamic formalism \cite{Pes97}.

For non-conformal repellers much less is known. The problem of existence of an ergodic measure with full dimension is solved for the class of Lalley-Gatzouras carpets and its nonlinear version \cite{LaGa92,Luz06,Luz10}. In \cite{GaPe97}, Gatzouras and Peres conjectured that such a measure is unique. However, in this note, we show that this may fail on linear Lalley-Gatzouras carpets. Such a phenomenon is known for some other examples of self-affine sets  constructed by K\"{a}enm\"{a}ki
and Vilppolainen \cite{KaVi10}.

To construct our example, let  $(X,\sigma_X)$ and $(Y,\sigma_Y)$ be  one-sided full shifts  over finite alphabets $\A$ and $\B$, respectively. Let $\pi:X\to Y$ be a $1$-block factor map, i.e., there is a map $\widetilde \pi:\A\to \B$  such that $$
\pi(x)=(\widetilde\pi(x_i))_{i=1}^\infty,\quad x=(x_i)_{i=1}^\infty\in X.$$
    Let $\phi: X\to \R$ and $\psi: Y\to \R$ be two positive functions which are constants over the cylinders  of first generation of $X$ and $Y$ respectively, i.e.,
 $$
 \phi(x)=\phi(x_1),\quad \psi(y)=\psi(y_1)
 $$
  for each $x=(x_i)_{i=1}^\infty\in X$ and  $y=(y_i)_{i=1}^\infty\in Y$.
 Furthermore, assume that $\phi(x)>\psi(\pi(x))$ for all $x\in X$.

 Define
\begin{equation}
\label{e-1}
P(\phi,\psi)=\sup
\left\{
\frac{h_\mu(\sigma_X)- h_{\mu\circ \pi^{-1}}(\sigma_Y)}{\int\phi \, d\mu}+
\frac{h_{\mu\circ \pi^{-1}}(\sigma_Y)}{\int\psi\circ \pi  \,d\mu}
\right\},
\end{equation}
where the supremum is taken over the collection $M(X,\sigma_X)$ of all $\sigma_X$-invariant Borel probability measures  on $X$. Here $h_\mu(\sigma_X)$ stands for the measure-theoretic entropy of $\sigma_X$ with respect to $\mu$ (cf. \cite{Pes97, Wal82}). Since the entropy maps $\mu\mapsto h_\mu(\sigma_X)$ and $\mu\mapsto h_{\mu\circ \pi^{-1}}(\sigma_Y)$ are upper semi-continuous, the supremum is attained on $M(X,\sigma_X)$. Moreover, since $\phi(x)$ and $\psi(y)$ only depend on the first coordinates of $x$ and $y$, the supermum can be only attained at Bernoulli measures in $M(X,\sigma_X)$ \footnote{As a related result,   Luzia \cite{Luz10} proved recently that  the supremum in (\ref{e-1}) always can be attained at ergodic measures when $\phi$ and $\psi$ are assumed to be  general positive H\"{o}lder continuous functions.
}.

In the next section, we construct an example to show that  in the above general setting, there may have two different Bernoulli measures in  $M(X,\sigma_X)$  attaining the supermum in \eqref{e-1}, which leads to a counter-example to Gatzouras and Peres conjecture on Lalley-Gatzouras carpets (see Section 3).

In fact, Gatzouras and Peres raised the wider conjecture claiming that if  $f$ is a smooth expanding map, then any compact invariant set $K$ which
satisfies specification carries a unique ergodic invariant measure $\mu$ of full dimension.
Moreover, $\mu$ is mixing for $f$. This conjecture was proved to be true in some special cases, e.g., as we said when $f$ is a conformal $C^{1+\alpha}$ map on smooth Riemanian manifolds \cite{GaPe97}, and also when $f$ is a linear diagonal endomorphism on the $d$-torus \cite{Fen11}. In particular, it is true for Bedford-McMullen self-affine carpets and sponges \cite{Bed84, McM84, KePe96} and some sofic self-affine sets \cite{Pet03,Ya09,Ol10}.

The same kind of questions have been studied on horseshoes. It is proved in \cite{MaMcC} that  for  nonlinear horseshoes there may be no ergodic measure of full dimension, while such a measure exists for linear horseshoes \cite{BaWo}, but may be not unique \cite{Rams}.

\section{An example}
Let $M(Y, \sigma_Y)$ denote the collection of all $\sigma_Y$-invariant Borel probability measures on $Y$. Notice that
$$
P(\phi,\psi)=\sup_{\nu\in \M(Y,\sigma_Y)} P(\phi,\psi,\nu),
$$
where
$$
P(\phi,\psi,\nu)=\frac{h_\nu(\sigma_Y)}{\int\psi \,d\nu}+P(\phi,\nu),\quad P(\phi,\nu)=\sup_{\substack{\mu\in \M(X,\sigma_X), \\\mu\circ\pi^{-1}=\nu}}\frac{h_\mu(\sigma_X)- h_\nu(\sigma_Y)}{\int\phi \, d\mu}
$$
for $\nu\in M(Y, \sigma_Y)$. Since $\phi(x)$ and $\psi(y)$ only depend on the first coordinates of $x$ and $y$,  $P(\phi,\psi,\nu)$ can  only be  maximized at  Bernoulli measures $\nu$ in $M(Y,\sigma_Y)$.

We make the following assumptions:

\begin{itemize}
\item[(1)] $\phi(x)\equiv \lambda>0$ on $X$ for some constant $\lambda$;

\item[(2)] $\B=\{a,b\}$, $\A=\{1,\cdots,\ell_a,\ell_a+1,\cdots,\ell_a+\ell_b\}$, $\widetilde\pi(\{1,\cdots,\ell_a\})=\{a\}$, $\widetilde\pi(\{\ell_a+1,\cdots,\ell_a+\ell_b\})=\{b\})$, where $\ell_a,\ell_b\in \N$.
\end{itemize}

Then, since due to our assumption we have $\int\phi \, d\mu=\lambda$ for all $\mu\in M(X,\sigma_X)$, the Ledrappier-Walters relativized variational principal \cite{LeWa77} yields
$$
 P(\phi,\nu)=\frac{1}{\lambda}(\log (\ell_a)\nu([a])+\log(\ell_b)\nu([b])),
$$
where $[c]:=\{y=(y_i)_{i=1}^\infty\in Y:\; y_1=c\}$ for $c\in \B$. Setting $x=\nu([a])$ and $H(x)=-x\log(x)-(1-x)\log(1-x)$, we thus have for all Bernoulli measures $\nu\in M(Y,\sigma_Y)$,
\begin{equation}
\label{e-a}
P(\phi,\psi,\nu)=f(x)=\frac{1}{\lambda}(\log (\ell_a/\ell_b)x+\log(\ell_b)) +\frac{H(x)}{\psi_ax+ \psi_b(1-x)},
\end{equation}
where $\psi_a$ and $\psi_b$ stand for the constant values of $\psi$ over $[a]$ and $[b]$ respectively.

A counter-example will appear if we find $\lambda, \ \ell_a, \ \ell_b, \psi_a$, and $\psi_b$ such that  $f$ attains its maximum for at least two values of $x$ in $[0,1]$.

Setting $U=\displaystyle\frac{\psi_b}{\lambda}\log (\ell_a/\ell_b)$ and $V=\displaystyle\frac{ \psi_a-\psi_b}{\psi_b}$, the problem transfers to finding $U\in \R$, $V\in(-1,\infty)$ and $M\geq 0$ such that
$$
g(x)=Ux-M+\frac{H(x)}{1+Vx}\le 0,\quad \forall\ x\in [0,1]
$$
and $g(x)=0$ has more than one solution  in $[0,1]$. We can seek for a quadratic polynomial $F(x)=A-B(x-1/2)^2$ with $A, B>0$ such that
\begin{itemize}
\item[(i)]
$F(x)\ge H(x)$ for all $x\in [0,1]$; and
\item[(ii)] the equation $F(x)=H(x)$ has more than one solution  in $[0,1]$.
 \end{itemize}
 Due to the common symmetry properties of $F$ and $H$ with respect to $x=1/2$ and the concavity of these functions, this will be the case if we make sure that the curvature of $F$ at $1/2$ is larger than that of $H$ at $1/2$ and $\inf_{x\in [0,1]}(F(x)-H(x))=0$. Recalling that the curvature of a smooth function $h(x)$ being given by $$\mathcal K_h(x)=\displaystyle \frac{|h''(x)|}{(1+(h'(x))^2)^{3/2}},$$
 we have $\mathcal K_H(1/2)=4$ and $\mathcal K_F(1/2)=2B$. Thus we get the following necessary and sufficient condition to guarantee that (i)-(ii) hold:
 \begin{equation}
 \label{e-l}
  B>2,\quad  A=\max_{0\leq x\leq 1}(B(x-1/2)^2+H(x)).
\end{equation}

Now take a pair of  numbers $A,B$ so that  \eqref{e-l} holds.  Then the identity
$$
-(Ux-M)(1+Vx)=A-B(x-1/2)^2
$$
yields
$$
\begin{cases}
UV=B,\\
MV-U=B,\\
M=A-B/4.
\end{cases}
$$
 This forces
$$
(A-B/4)V^2-BV-B=0.
$$
The positive root of the above equation is
\begin{equation}
\label{e-m}
V=\displaystyle \frac{2B+4\sqrt{AB}}{4A-B}.
\end{equation}
Then, using the equality $UV=B$ yields
\begin{equation}
\label{e-n}
U=\sqrt{AB}-B/2.
\end{equation}
Next take
\begin{equation}
\label{e-o}
\psi_b=1,\ \psi_a-\psi_b=V,
\end{equation}
and take positive integers  $\ell_a$, $\ell_b$
such that
\begin{equation}
\label{e-q}
\log(\ell_a/\ell_b)>\frac{1+V}{V} B.
\end{equation}
In the end,  take $\lambda$ such that
\begin{equation}
\label{e-p}
\frac{\log(\ell_a/\ell_b)}{\lambda}=U,\ i.e.,\  \lambda=\frac{ \log(\ell_a/\ell_b)}{U}=\log(\ell_a/\ell_b)\frac{V}{B}.
\end{equation}
According to \eqref{e-q}-\eqref{e-p}, $\lambda>1+V$ and thus  $\phi\equiv\lambda>\max(\psi)=\max (\psi_a,\psi_b)$.

Then for the above constructed  $\lambda, \ \ell_a, \ \ell_b, \psi_a$, and $\psi_b$, the function $f(x)$ defined in \eqref{e-a} attains its supremum at two different points $x$ in $[0,1]$. This yields an example that the supermum in \eqref{e-1} is attained at two different Bernoulli measures in $M(X,\sigma_X)$.

In the end, we provide a more concrete example for $\lambda, \ \ell_a, \ \ell_b, \psi_a$, and $\psi_b$.

\begin{ex}
\label{ex-1}
Set $$
B=3\log 2\approx 2.07944$$ and $$A=\log 3-\frac{7}{12}\log 2\approx 0.69427643.$$
One can check that \eqref{e-1} holds for such $A$ and $B$. Indeed, the supermum in defining $A$ is attained at $x=1/3$.
Then
$$U=\sqrt{AB}-B/2\approx 0.16182292, \quad V=\displaystyle \frac{2B+4\sqrt{AB}}{4A-B}\approx 12.8501046.$$
Take $$\psi_a=1+V\approx 13.8501046, \quad \psi_b=1$$
and
$$\ell_a=150,\quad \ell_b=1,\quad \lambda=\log (\ell_a/\ell_b)\cdot \frac{V}{B}\approx 30.9636922.$$
\end{ex}
\section{Application to Gatzouras-Lalley carpets}
\label{S-3}

Let $\lambda, \ \ell_a, \ \ell_b, \psi_a$, and $\psi_b$ be constructed as in Example \ref{ex-1}. Notice that
$$3\exp(-\lambda)\ell_a<3\cdot e^{-30}\cdot 150< 1,\quad \exp(-\psi_a)+\exp(-\psi_b)<2e^{-1}<1.$$
Then we can build a Gatzouras-Lalley carpet in the unit square as the attractor $K$ of the IFS $\{S_{a,r}:1\le r\le \ell_a\}\cup\{S_{b,s}:1\le s\le \ell_b\}$, where
$$
\begin{cases}
S_{a,r}(x,y)=(\exp(-\lambda)x, \exp(-\psi_a) y)+ (2 r\exp(-\lambda),0),\ 1\le r\le \ell_a,\\
S_{b,s}(x,y)=(\exp(-\lambda)x, \exp(-\psi_b) y)+ (2 r\exp(-\lambda),1-\exp(-\psi_b)),\ 1\le s\le \ell_b.
\end{cases}
$$
Gatzouras and Lalley \cite{LaGa92} proved  that  the Hausdorff dimension of $K$ is equal to $P(\phi,\psi)$, which is attained by some  Bernoulli measure on $X$.   The previous section shows that such a measure is not unique; in our example there are exactly two such measures.

\noindent{\bf Acknowledgements}. The authors would like to thank Antti K\"{a}enm\"{a}ki for  pointing out  the references \cite{BaWo, KaVi10,Rams} and some related results after they submitted the first version of this paper. Both authors were partially supported by the   France/Hong Kong joint research scheme {\it PROCORE} (projects 20650VJ,  F-HK08/08T). Feng was also partially
supported by the RGC grant and the Focused Investments Scheme in the CUHK.

\end{document}